\documentclass[12pt]{article}

\usepackage[T2A]{fontenc}
\usepackage[utf8]{inputenc}

\usepackage{amsfonts}
\usepackage{amssymb}
\usepackage{amsmath}
\usepackage{amsthm}

\def \le {\leqslant}
\def \ge {\geqslant}

\theoremstyle{plain}
 
\begin{document}
\begin{Huge}
 \centerline{\"{U}ber die Winkel zwischen Unterr\"{a}umen}
\end{Huge}
\vskip+0.5cm
\begin{Large}
\centerline{Nikolay Moshchevitin\footnote{Diese Arbeit wurde an der Nationalen Pazifik-Universit\"{a}t durchgef\"{u}hrt und durch RNF unterstützt (Grant no. 18-41-05001)  }
}
\end{Large}
 
\vskip+1cm
\begin{small}
 {\bf Abstract.}
 We prove  a metric statement about approximation of  a $n$-dimensional   linear subspace $A$ in $\mathbb{R}^d$ by $n$-dimensional rational subspaces.  
  We consider the problem of finding  a rational subspace $B$  of bounded height $H=H(B)$ for which the angle of inclination
$\psi (A,B) $ is small in terms of $H$.
 In the simplest case $d=4, n=2$ we give a partial solution of a problem formulated by W.M. Schmidt in 1967.
 
 \vskip+0.3cm
 {\bf AMS Subject Classification:} 11J13.
 
  \vskip+0.3cm
 {\bf Keywords:} Diophantine Approximation, subspaces, angles of inclination.
 \end{small}
 \vskip+2cm
 
 Das Ziel dieser Arbeit ist  es,
 einen neuen metrischen Satz 
 \"{u}ber   Approximationen  der $n$-dimensionalen reelen Unterr\"{a}ume 
 von $\mathbb{R}^d$
 durch 
 $n$-dimensionale
 rationale Unterr\"{a}ume 
 zu beweisen.
 Die wichtigsten Fragen zur Theorie dieser Approximationen wurden im Jahr 1967 von W.M. Schmidt formuliert \cite{s1}.
 Hier beschränken wir uns auf den einfachsten Fall
 $d=4, n=2$, den Schmidt betrachtet hat.
 Wir formulieren  unseren Hauptsatz (Satz 1)
 in Sektion 1 und in  den Sektionen 2 bis 6 geben wir den vollst\"{a}ndigen Beweis.
 In Sektion 7 formulieren wir einen allgemeinen Satz ohne Beweis.

  \vskip+0.3cm
 {\bf 1. Winkel zwischen Unterr\"{a}umen.}
 \vskip+0.3cm
 
 Wir betrachten den $d$-dimensionalen euklidischen Raum $\mathbb{R}^d$ 
 mit Koordinaten $ (z_1,...,z_d)$ versehen
 mit dem Skalarprodukt
 $$
 \langle \pmb{w},\pmb{z}\rangle = w_1z_1+...+w_dz_d.
 $$
 Weiters
 bezeichnen  wir
 die euklidische Norm des Vektors $\pmb{z}$ mit
 $$
 |\pmb{z}| = 
 \sqrt{z_1^2+...+z_d^2}.
 $$
 Sei 
 $$
 \rho ( \mathcal{A}_1, \mathcal{A}_2) = 
 \min_{\pmb{w}\in \mathcal{A}_1,\, \pmb{z}\in \mathcal{A}_2}
 \,\, |\pmb{w} - \pmb{z}|
 $$
 der
 euklidische Abstand
 zwischen  den Mengen $\mathcal{A}_1$ und $\mathcal{A}_2$.  Mit
 $$
 \mathcal{U}_\varepsilon (\mathcal{A})= 
 \{
 \pmb{z}\in \mathbb{R}^d:\,\,\,
 \rho (\pmb{z}, \mathcal{A})\le \varepsilon\}
 $$
 bezeichnen wir die $\varepsilon$-Umgebung der Menge $\mathcal{A} \subset \mathbb{R}^d$.
 
 Gegeben sei das Gitter $\mathbb{Z}^d\subset \mathbb{R}^d$. 
 F\"{u}r jeden linearen $n$-dimensionalen Unterraum $ L\subset \mathbb{R}^d$
 betrachten wir das  Gitter $ \Gamma_L = L\cap \mathbb{Z}^d$.
 Der Unterraum $L$ hei\ss t {\it rational}, wenn
 der Rang des Gitters $\Gamma_L $ 
 gleich der Dimension des Unterraums $L$ ist.
 F\"{u}r einen $n$-dimensionalen  rationalen Unterraum $L$
 definieren wir die
 Höhe  $H(L)$ des Unterraums $L$ 
 als die Gitterdiskriminante des Gitters $ \Gamma_L$.
 Es ist klar, dass $H(L)\ge 1$ und $ (H(L))^2\in \mathbb{Z}$ gelten muss.
 Es sei
 $$
 \hbox{\got S} =
 \{ \pmb{z} = (z_1,...,z_d) \in \mathbb{R}^d:\,\,\,
 z_1^2+...+z_d^2 =1\}
 $$
 die Oberfl\"{a}che  der Einheitskugel im $\mathbb{R}^d$. Wir definieren den Winkel 
$\psi ( L_1, L_2)$
 zwischen zwei 
 Unterr\"{a}umen $L_1$ und $L_2$ durch
 $$
 \psi ( L_1, L_2)  = \min_{\pmb{w}\in \omega_{L_1},\, \pmb{z}\in \omega_{L_2}}
 \,\, \sigma (\pmb{w}, \pmb{z}),
$$
 wobei $
 \omega_L = L\cap \hbox{\got S},
 $
 und 
 $
 \sigma (\pmb{w}, \pmb{z})$  den Winkel
 zwischen den Vektoren 
 $\pmb{w}  $ und $\pmb{z}$ bezeichnet, also
 $$
 \cos \left(
 \sigma (\pmb{w}, \pmb{z})
 \right) =
 \frac{\langle\pmb{w},\pmb{z}\rangle}{|\pmb{w}|\cdot|\pmb{z}|}.
 $$
 Daher entspricht $  \sigma (\pmb{w}, \pmb{z})$  dem Abstand 
zwischen  den Punkten $\pmb{w}$ und $\pmb{z}$
 auf  
 $
 \hbox{\got S}$ 
(in der Notation der Arbeit \cite{s1} haben wir daher
$ \psi ( L_1, L_2) =  \psi_1 ( L_1, L_2)$). 
 Es ist klar, dass f\"{u}r  $\pmb{w}, \pmb{z}\in \hbox{\got S} $  die Ungleichung 
 $$
 |\pmb{w} - \pmb{z}| \le   \sigma (\pmb{w}, \pmb{z})
 $$
 gilt.

1967 bewies W. M. Schmidt \cite{s1} viele Resultate \"{u}ber
 diophantische Approximationen mit Unterr\"{a}umen, darunter die folgenden
 S\"{a}tze A und B bewiesen
 (siehe Theorem 12 und  Theorem 16 sowie Corollary 2 aus \cite{s1}).
 
  \vskip+0.3cm
  {\bf Satz A.}\,{\it Seien $n,k$  positive ganze Zahlen mit $ n+k \le d$. 
  Es gibt eine positive Konstante $C_1$ mit der folgenden Eigenschaft.
  
  F\"{u}r jeden $n$-dimensionalen Unterraum $A$ und jedes $H>1$ 
  gibt es  einen rationalen $k$-dimensionalen Unterraum $B$ mit
    \vskip+0.3cm
    
    {\bf (i)} \,\,\,\, $H(B) \le H$,

    {\bf (ii)} \,\,\,  $\psi (A,B) \le C_1 H^{-\frac{(d-1)n}{(d-n)(d-k)}}\cdot (H(B))^{-\frac{d-1}{d-k}}.$
    
        \vskip+0.3cm
    
    Insbesondere wenn $A\cap B'=\{\pmb{0}\}$
    f\"{u}r alle $k$-dimensionalen rationalen Unterr\"{a}ume $B'$ gilt,
    gibt es unendlich viele $k$-dimensionale rationale Unterr\"{a}ume $B$ mit
    $$
    \psi (A,B) \le C_1 
     (H(B))^{-\frac{d(d-1)}{(d-n)(d-k)}}.
    $$ 
    }
    
      \vskip+0.3cm
  
    {\bf Satz B.}\,{\it Seien $n,k$  positive ganze Zahlen mit $ n+k \le d$.
    Dann gibt es ein positives $C_2$ und einen $n$-dimensionalen 
    Unterraum  $A$ mit der folgenden Eigenschaft. F\"{u}r jeden rationalen 
    $k$-dimensionalen Unterraum $B$ hat man
    $$
    \psi (A,B) \ge C_2  (H(B))^{-\frac{k(d-k)+1}{d+1-n-k}}
    .$$
    }
 
 Im einfachsten Fall $k=1$  ist Satz 1 analog zum 
 Schubfachprinzip von Dirichlet (siehe \cite{c}, Ch.1, Theorem VI, 
 oder \cite{s}, Ch.2). Ebenfalls f\"{u}r $ k =1$ ist Satz B   analog zu der Behauptung aus dem Buch \cite{c} 
 (siehe \cite{c}, Ch.1, Theorem VIII), die     auf Perron \cite{pe}  zur\"{u}ck geht.
 Der  Beweis von Satz A  f\"{u}r $ k > 1$  basiert auf mehreren Anwendungen des 
 Minkowskischen Gitterpunktsatzes.
 Der Unterraum $A$ aus Satz B wurde mittels algebraischer Zahlen konstruiert.
 
 Betrachten wir  nun den  Fall $d= 4, n=k=2$. Dann 
 lauten die Aussagen der S\"{a}tze A und  B kurz wie folgt:
 
 a) 
 f\"{u}r jeden $2$-dimensionalen  Unterraum $A$ und jedes $H>1$ 
  gibt es  einen rationalen   $2$-dimensionalen  Unterraum  $B$ mit
$$
 \psi (A,B) \le C_1 H^{- 3}
 \,\,\,\,\,\,
 \text{und}
 \,\,\,\,\,\,
 H(B) \le H.
 $$
 
 b) es gibt
 einen $2$-dimensionalen  Unterraum 
  $A$ mit 
    $$
    \psi (A,B) \ge C_2  (H(B))^{-5}
    \,\,\,\,\,
    \text{f\"{u}r alle $2$-dimensionalen rationalen Unterr\"{a}ume}\,\,\,
    B.
    $$

    In \cite{s1},\S 16 stellt Schmidt die
  Frage nach dem bestm\"{o}glichen  Exponenten in  den Behauptungen  a) und b).
    Nun formulieren wir eine Versch\"{a}rfung der Behauptung b).

      \vskip+0.3cm
  
       {\bf Satz 1.}\,\,{\it   
    Sei
    $\omega (j) \ge 0, j \ge 1,$ 
    eine positive  monoton fallende Funktion
 mit
     \begin{equation}\label{reiche0}
 \sum_{j=1}^\infty\,  j \, \omega \left(\sqrt{j}\right)<+\infty. 
 \end{equation}

    Dann  gibt es f\"{u}r fast jeden 2-dimensionalen Unterraum $A
    \subset\mathbb{R}^4$  
    eine positive Konstante $C(A)$ mit
    $$
     \psi (A,B) \ge C(A) \, \omega (H(B)) \,\,\,\,\,
    \text{f\"{u}r alle $2$-dimensionalen rationalen Unterr\"{a}ume}\,\,\,
    B
    \,\,\,
    \text{in}
    \,\,\,
    \mathbb{R}^4.
    $$
     Insbesondere existiert  f\"{u}r jedes $\varepsilon >0$
     f\"{u}r fast jeden  2-dimensionalen Unterraum $A\subset \mathbb{R}^4$ 
    eine positive Konstante $C(A,\varepsilon)$
    sodass
    $$
     \psi (A,B) \ge C(A,\varepsilon) \,(H(B))^{-4-\varepsilon}
     \,\,\,\,\,
    \text{f\"{u}r alle $2$-dimensionalen rationalen Unterr\"{a}ume}\,\,\,
    B
     \,\,\,
    \text{in}
    \,\,\,
    \mathbb{R}^4\,\,\,\text{gilt}.
    $$
    }
    Auf der Menge der Unterr\"{a}me
   existiert ein  invariantes Ma\ss \,  
(siehe zum Beispiel \cite{sw}, Ch. 13), 
dessen Verwendung kann hier allerdings umgangen werden.
In der n\"{a}chsten Sektion  erkl\"{a}ren wir 
die Bedeutung von "für fast alle Unterr\"{a}ume". Das $l$-dimensionale Lebesgumaß
in $\mathbb{R}^l$  sei nach folgend mit $\mu_l(\cdot)$ bezeichnet.

  \vskip+0.3cm
{\bf 2.  Projektionen und die Mengen von Unterr\"{a}umen.}
  \vskip+0.3cm

Sei $L$ ein zweidimensionaler linearer Unterraum des Raumes $\mathbb{R}^4$ und sei
\begin{equation}\label{baas}
\pmb{q}_1= (q_{1,1},q_{1,2},q_{1,3},q_{1,4}),\,\,\,\,\,
\pmb{q}_2= (q_{2,1},q_{2,2},q_{2,3},q_{2,4})
\end{equation}
eine Basis des Raums $L$.
Wir definieren 
$$
p_{i,j} =
\left|
\begin{array}{cc}
q_{1,i}& q_{1,j}
\cr
q_{2,i}& q_{2,j}
\end{array}
\right|,\,\,\, i,j \in \{ 1,2,3,4\}
$$
und
betrachten  die  homogenen Plücker-Koordinaten 
$$
\hbox{\got p} =(
p_{1,2} : p_{1,3} : p_{1,4} : p_{2,3} : p_{2,4} : p_{3,4})
$$
des Unterraums $L$.
Die Menge der zweidimensionalen Unterräume des $\mathbb{R}^4$
 kann mit  der
 Graßmann-Mannigfaltigkeit $ \hbox{\got G}=\hbox{\got G}(2,2)$
  identifiziert werden (siehe \cite{hp}).
  Die  Graßmann-Mannigfaltigkeit $ \hbox{\got G}$
ist eine Mannigfaltigkeit im reellen projektiven Raum $\mathbb{P}^5$, die   mittels Plücker-Koordinaten 
definiert werden kann:
$$
\hbox{\got G} = \{ \hbox{\got p} \in \mathbb{P}^5:\,\,\,\,
p_{1,2}p_{3,4}-p_{1,3}p_{2,4}+p_{1,4}p_{2,3}=0\}.
$$ 
Diese Darstellung liefert eine Bijektion zwischen der Menge der zweidimensionalen Unterräume und 
$
\hbox{\got G} .
$ 
In der Folge 
identifizieren wir die Unterr\"{a}ume von $\mathbb{R}^4$ mit den Punkten von $\hbox{\got G}$. 
Weiters 
betrachten wir  die Menge
$$
\hbox{\got G}^{[0]} =
\{ L \in 
\hbox{\got G}:\,\,\,
\exists\, i,j \,\,\,\text{mit}\,\,\, p_{i,j} = 0\}.
$$

Wir schreiben Punkte in 
$\mathbb{R}^6$ als
 $$ \pmb{p} = ( p_{1,2}, p_{1,3}, p_{1,4}, p_{2,3}, p_{2,4}, p_{3,4})$$
und  betrachten f\"{u}r
$i< j$
 die affinen Unterr\"{a}ume
$$
\hbox{\got E}_{i,j} =\{ \pmb{p} =( p_{1,2}, p_{1,3}, p_{1,4}, p_{2,3}, p_{2,4}, p_{3,4})\in \mathbb{R}^6:\,\,\,\,
p_{i,j} = 1\} \subset \mathbb{R}^6
$$
sowie die  Mengen
$$
\hbox{\got G}_{i,j} =
\{ \pmb{p} \in 
\hbox{\got E}_{i,j}
:\,\,\,\,
p_{1,2}p_{3,4}-p_{1,3}p_{2,4}+p_{1,4}p_{2,3}=0\}
\subset \hbox{\got E}_{i,j},
$$
$$
\hbox{\got V}_{i,j} =\{ \pmb{p} \in  \hbox{\got E}_{i,j} :\,\,\, |p_{l,k}| \le 1,\,\,\, l<k\},\,\,\,\,\,\,\,
\hbox{\got W}_{i,j} =
\{ \pmb{p} \in  \hbox{\got V}_{i,j} :\,\,\, p_{1,2}p_{3,4}-p_{1,3}p_{2,4}+p_{1,4}p_{2,3}=0\}
$$
 und
$$
\hbox{\got W} =\bigcup_{i<j}
\hbox{\got W}_{i,j} .
$$
Dann gibt es eine Karte
$$
\pi:\,\,
\hbox{\got G} \rightarrow
\hbox{\got W}    ,\,\,\,\,\,\,\,\,\,\,  L\,\, \text{mit Koordinaten}\,\,\hbox{\got p} \,\,\,\,\,\mapsto\,\,\,\,\, \pmb{p},
$$
die 
im Allgemeinen nicht
bijektiv ist,
sondern lediglich
die eingeschränkte Karte
$$
\pi
:\,\,
\hbox{\got G} 
\setminus \left(\pi^{-1}\left( \partial \,\hbox{\got W} \right)\right)
\rightarrow
\hbox{\got W} \setminus  \left( \partial \,\hbox{\got W} \right)
$$ 
ist bijektiv. Hier ist
$$
\partial \,\hbox{\got W} =  \hbox{\got W} \cap \left(
\bigcup_{i<j}\,\,
\{ \pmb{p} \in  \hbox{\got V}_{i,j} :\,\,
\exists (k,l)\,\,\,\text{mit}\,\,\, p_{k,l}= -1\}
\right)
$$
die "Grenze" der Menge $\hbox{\got W} $.

Sei  
$i_1,i_2,i_3,i_4$ eine Permutation der Folge $1,2,3,4$ mit $ i_1<i_2, i_3<i_4$ und
 $$\pmb{p} =(p_{1,2}, p_{1,3}, p_{1,4}, p_{2,3}, p_{2,4}, p_{3,4})\in \hbox{\got W}_{i_1,i_2} .$$
Dann ist
$$
p_{i_1,i_2}p_{i_3,i_4}-p_{i_1,i_3}p_{i_2,i_4}+p_{i_1,i_4}p_{i_2,i_3}=0,
$$
und  f\"{u}r $p_{i,j}$ 
hat man
\begin{equation}\label{ab}
p_{i_1,i_2} =  1,\,\,\,
p_{i_1,i_3} =\ell_{1,2},\,\,\,
 p_{i_1,i_4}=\ell_{2,2} ,  \,\,\,
 p_{i_2,i_3}= -\ell_{1,1} ,\,\,\,
  p_{i_2,i_4}= -\ell_{2,1}, \,\,\,
  p_{i_3,i_4} = \ell_{1,1}\ell_{2,2} -\ell_{2,1} \ell_{1,2}  
\end{equation}
mit reellen $ \ell_{i,j} \in [-1,1]$ und $ \Delta (\pmb{\ell}) =\ell_{1,1}\ell_{2,2} -\ell_{2,1} \ell_{1,2}\in [-1,1]$. 

 Nun möchten wir eine   einfache Behauptung formulieren, die sofort aus den Eigenschaften der Karte 
$ \pi$ folgt.
 
   \vskip+0.3cm
 {\bf Hilfssatz 1.}\,\,{\it
Sei $ L\subset \mathbb{R}^4$ ein $2$-dimensionaler Unterraum. Dann gibt es eine
Permutation $i_1,i_2, i_3, i_4$, so dass man in  den Koordinaten
$$
x_1 = z_{i_1},\,\,\,\,\, x_2 = z_{i_2},\,\,\,\,\, y_1 = z_{i_3},\,\,\,\,\, y_2 = z_{i_4}
$$
die Identit\"{a}t
\begin{equation}\label{ell}
L =\{\pmb{z}= (x_1,x_2,y_1,y_2):\,\,\, y_j = \ell_{j,1} x_1 +
\ell_{j,2} x_2,\,\,\,   j =1,2\}
\end{equation}
mit $\ell_{i,j}$ aus (\ref{ab}) hat.
Insbesondere ist
\begin{equation}
\label{le}
\max_{1\le i,j\le 2} 
|\ell_{i,j}|
\le 1.
\end{equation}
}
  \vskip+0.3cm

Beweis. Wir betrachten  nur den Fall
$\pmb{p} \in \hbox{\got W}_{1,2}$ und  $(i_1,i_2,i_3,i_4) = (1,2,3,4)$.
Dann sind
$$
w_1^{[1]} = \frac{
p_{3,2}
}{p_{1,2}} = -p_{2,3} = \ell_{1,1},\,\,\,\,\,
w_2^{[1]} = \frac{
p_{1,3}
}{p_{1,2}} = p_{1,3} = \ell_{1,2},\,\,\,\,\,
w_3^{[1]} =  
-1,\,\,\,\,\,
w_4^{[1]} = 0
$$
und
$$
w_1^{[2]} = \frac{
p_{4,2}
}{p_{1,2}} = -p_{2,4} = \ell_{2,1},\,\,\,\,\,
w_2^{[2]} = \frac{
p_{1,4}
}{p_{1,2}} = p_{1,4} = \ell_{2,2},\,\,\,\,\,
w_3^{[2]} =  
0,\,\,\,\,\,
w_4^{[2]} = -1
$$
 L\"{o}sungen 
des Systems
$$
\begin{cases}
q_{1,1}w_1+q_{1,2}w_2+q_{1,3}w_3 + q_{1,4}w_4 = 0\cr
q_{2,1}w_1+q_{2,2}w_2+q_{2,3}w_3 + q_{2,4}w_4 = 0
\end{cases},
$$
und die Vektoren (\ref{baas})  
erf\"{u}llen die Gleichungen aus
(\ref{ell}).$\Box$

  \vskip+0.3cm

Wir definieren nun die Mengen
 \begin{equation}\label{eee}
 \mathcal{E} =
\{ \pmb{\ell} =(\ell_{1,1},\ell_{1,2},\ell_{2,1},\ell_{2,2}) \in \mathbb{R}^4:\,\,\,  \ell_{i,j}   \in [-1,1],\,\,\,| \Delta (\pmb{\ell}) |\le 1\} \subset\mathbb{R}^4,
\end{equation}
$$
\overline{\mathcal{E} } = \{ (i,j), \,\,i<j,\,\, i,j \in \{1,2,3,4\}\}\times \mathcal{E}
\subset \{ (i,j), \,\,i<j,\,\, i,j \in \{1,2,3,4\}\} \times\mathbb{R}^4
,
$$
und die Karten
$$
\nu_{i,j}:\,\, 
\hbox{\got G}_{i,j}
\,\, \rightarrow \,\,\mathbb{R}^4,\,\,\,\,\,
 \pmb{p}  \,\,\mapsto\,\,  \pmb{\ell};
 \,\,\,\,\,\,
 \nu_{i,j}(
\hbox{\got W}_{i,j}) = {\mathcal{E} },
 $$
 und
 $$
\nu:\,\,
  \pmb{p} \in \hbox{\got W}_{i,j}\setminus \left( \partial \,\hbox{\got W}_{i,j}\right)\,\,\mapsto\,\,  ((i,j), \pmb{\ell})\in  \overline{\mathcal{E} }, 
$$
wobei die reellen Zahlen  $\ell_{i,j}$   durch (\ref{ab}) definiert wurden und
$$
\partial \,\hbox{\got W}_{i,j} =
\hbox{\got W}_{i,j}\cap \{ \pmb{p} \in  \hbox{\got V}_{i,j} :\,\,\,
\exists (k,l)\neq (i,j)\,\,\,\,\text{mit}\,\,\, |p_{k,l}|= 1\}
$$
  die "Grenze" der Menge $\hbox{\got W}_{i,j} $ ist.
Betrachten wir die Mengen
$$
\hbox{\got G}^{[1]}
=
\hbox{\got G} \setminus  \left(\hbox{\got G}^{[0]}   \cup
\left(\pi^{-1}\left( \bigcup_{i<j}
\partial \,\hbox{\got W}_{i,j} \right)\right) \right)  
$$
und
$$
 \overline{\mathcal{E}}^{[1]} =\{ 
 ((i,j),\pmb{\ell})\in  \overline{\mathcal{E} }:\,\,\,\,
 0< |\ell_{i,j}|,  < 1,\,\,\,i,j = 1,2,\,\,\,
 0< |
 \Delta (\pmb{\ell}) |<1\},
$$
so haben wir die bijektive Karte
$$
  \nu\circ \pi:\,\,\,
\hbox{\got G}^{[1]} 
  \rightarrow
\overline{\mathcal{E} }^{[1]}
,\,\,\,\,\,\,\,\,\,\,\,\, (\, L\,\, \text{mit Koordinaten}\,\,\hbox{\got p}\,) \,\mapsto\,((i,j),\pmb{\ell})
$$
definiert, 
wobei $
\pi (L) \in \hbox{\got E}_{i,j}$ und $\ell_{i,j}$   durch (\ref{ab}) definiert wurden.

F\"{u}r eine Menge $ \hbox{\got A} \subset \hbox{\got G}$ definieren wir die Projektionen
$$
\hbox{\got A}_{i,j} = \pi ( \hbox{\got A} )\cap \hbox{\got W}_{i,j} 
\,\,\,\,\,\text{und}
\,\,\,\,\,
\mathcal{A}_{i,j} = \nu_{i,j} \left(  \hbox{\got A}_{i,j} \right).
$$
Nun erkl\"{a}ren wir die Bedeutung   des Satzes 1.
Die Behauptung gilt  "für fast alle Unterr\"{a}ume" wenn 
man
f\"{u}r die
Ausnahmemenge  $ \hbox{\got A} \subset \hbox{\got G}$    
$$
\mu_4  \left(  \mathcal{ A}_{i,j} \right) = 0\,\,\,\,\,
\text{f\"{u}r alle}
\,\,\,\,\,
i<j
$$
hat.
Wir
bemerken an dieser stelle noch, dass
\begin{equation}\label{00}
\mu_4 \left(
\nu_{i,j} \left(\partial \,\hbox{\got W}_{i,j} \right)\right) = \mu_4 \left(
\nu_{i,j}\left(\pi\left(\ \hbox{\got G}^{[0]}   \right) \cap \hbox{\got W}_{i,j}  \right)\right) =0,
\,\,\,\,\,\,\,
\forall\, i,j
\end{equation}
gilt.
Insbesondere ist die Grenze eine Nullmenge.

Wir brauchen mehrere Projektionen.
Wir definieren   die Karte
$$
\lambda_{i,j;k,l} :\, \hbox{\got E}_{i,j}  \rightarrow  \hbox{\got E}_{k,l},\,\,\,\,\,
(\,\pmb{p} = (p_{1,2} , p_{1,3},..., p_{3,4})
\,\, \text{mit}\,\, p_{i,j} = 1 \,)
\,
\mapsto\,(\,
\pmb{p}' = (p_{1,2}' ,p_{1,3}'..., p_{3,4}')
\,\,\text{mit}\,\, p_{k,l} = 1\,),
$$
$$\hbox{\got p} = 
 (p_{1,2} : p_{1,3} : ...: p_{3,4}) =
 \hbox{\got p}'=
 (p_{1,2}' : p_{1,3}' : ...: p_{3,4}')
,
$$
sodass   $\lambda_{i,j;k,l}$ die zentrale Projektion aus $\hbox{\got E}_{i,j}$ nach  $\hbox{\got E}_{k,l}$ ist . 
Weiters setzen wir
\begin{equation}\label{projektion}
\nu_{i,j;k,l} = \nu_{k,l}\circ \lambda_{i,j;k,l} \circ \nu_{i,j}^{-1}
:
\,\,\mathbb{R}^4 \rightarrow \mathbb{R}^4.
\end{equation}
Sei $ \delta >0$.
Wir definieren  die Mengen
 \begin{equation}\label{eee1}
 \mathcal{E} (\delta)=
\{ \pmb{\ell}\in \mathbb{R}^4:\,\,\, \delta\le|\ell_{i,j} |\le 1-\delta ,\,\,\, i,j = 1,2,\,\,\,\,\,\,\delta\le|\Delta (\pmb{\ell}) |\le 1-\delta\},
\end{equation}
$$
 \hbox{\got W}_{i,j} (\delta) =
 \nu_{i,j}^{-1} \left(\mathcal{E}(\delta)\right)=
 $$
  \begin{equation}\label{eee1q}
  =
 \{\pmb{p}:\,\,\,p_{i,j} = 1; \,\,\, \delta \le |p_{k,l}| \le 1 -\delta,\,\,\,  (k,l) \neq (i,j);\,\,\,\,\,
 p_{1,2}p_{3,4}-p_{1,3}p_{2,4}+p_{1,4}p_{2,3}=0\},
 \end{equation}
 $$
\overline{\mathcal{E} }  (\delta)= \{ (i,j), \,\,i<j,\,\, i,j \in \{1,2,3,4\}\}\times \mathcal{E} (\delta)
$$
und
 \begin{equation}\label{eee11}
\hbox{\got G} (\delta) = (\nu\circ \pi)^{-1} \left(\overline{\mathcal{E} }  (\delta)\right).
\end{equation}
Es ist nun klar, dass
 $\nu_{i,j;i,j} ({\mathcal{E} }  (\delta)) = 
 {\mathcal{E} }  (\delta) .
 $

 \vskip+0.3cm
 {\bf Hilfssatz 2.}\,\,{\it
  F\"{u}r $(i,j) \neq (k,l)$ hat man
 \begin{equation}\label{max}
   \max_{
   \pmb{\ell}'\in  \nu_{i,j;k,l} (\mathcal{E}  (\delta)) }\,\,
   \max_{1\le i,j\le 2}
|\ell_{i,j}'|
\le \frac{1}{\delta}.
\end{equation} }
 \vskip+0.3cm
 
 Beweis.
 Es gen\"{u}gt, die Karten
 $\nu_{1,2;k,l}$
  mit
 $(k,l) = (3,4)$ 
 oder
 $ (k,l) \in \{ (1,3), (1,4), (2,3), (2,4)\}$ 
 zu betrachten.
 
 Falls  $(k,l) = (3,4)$,  betrachten wir den
 Unterraum $ L = (\pi^{-1}\circ \nu_{1,2}^{-1}) (\ell)$
 mit $\ell \in \mathcal{E}  (\delta)$. Dann
 kann $L$  sowohl durch  die Gleichungen
 \begin{equation}\label{L1}
 \begin{cases}
 z_3 = \ell_{1,1}z_1+ \ell_{1,2}z_2,
\cr
 z_4 =\ell_{2,1}z_1+\ell_{2,2}z_2 ,
 \end{cases}\,\,\,\,\,\,\,\,\,\,
  \ell \in  \mathcal{E}  (\delta)
  \end{equation}
  als auch
    durch die anderen Gleichungen
  $$
 \begin{cases}
 z_1 = \ell_{1,1}'z_3+ \ell_{1,2}'z_4,
\cr
 z_2 =\ell_{2,1}'z_3+\ell_{2,2}'z_4 ,
 \end{cases}\,\,\,\,\,\,\,\,\,\,
  \ell' \in \nu_{1,2;3,4}( \mathcal{E}  (\delta))
  $$ definiert werden. Also ist
  $
  \left(
  \begin{array}{cc}
  \ell_{1,1}'&\ell_{1,2}'\cr
  \ell_{2,1}'&\ell_{2,2}'
  \end{array}
  \right)
  $
  die Inverse  der Martix von
  $
  \left(
  \begin{array}{cc}
  \ell_{1,1}&\ell_{1,2}\cr
  \ell_{2,1}'&\ell_{2,2}
  \end{array}
  \right)
  $ und
  $$
  \max_{i,j} |\ell_{i,j}'|
 \le
 \frac{\max_{i,j} |\ell_{i,j}|}{|\Delta(\ell)|}
 \le \frac{1-\delta}{\delta}< \frac{1}{\delta}.
 $$
 
 Die F\"{a}lle 
  $ (k,l) =(1,3), (1,4), (2,3), (2,4)$
  sind \"{a}hnlich. Wir betrachten den Fall   $ (k,l) =(1,3)$.
  Dann wird der Unterraum $L$ durch die Gleichungen (\ref{L1})
  definiert, und wieder
  k\"{o}nnen wir denselben
  Unterraum $L$ auch durch die Gleichungen 
  $$
   \begin{cases}
 z_2 = \ell_{1,1}'z_1+ \ell_{1,2}'z_3,
\cr
 z_4 =\ell_{2,1}'z_1+\ell_{2,2}'z_3 ,
 \end{cases}\,\,\,\,\,\,\,\,\,\,
  \ell' \in \nu_{1,2;1,3}( \mathcal{E}  (\delta))
  $$
    definieren. Daraus ergibt sich, dass
    $$
    \ell_{1,1}' 
   =
   -\frac{\ell_{1,1}}{\ell_{1,2}},\,\,\,\,\,\,
      \ell_{1,2}' 
   =
   \frac{1}{\ell_{1,2}},\,\,\,\,\,\,
      \ell_{2,1}' 
   =
   -\frac{\Delta(\ell )}{\ell_{1,2}},\,\,\,\,\,\,
   \ell_{2,2}' 
   =
   \frac{\ell_{2,2}}{\ell_{1,2}},
    $$
    und
    $ \max_{i,j}|\ell_{i,j}'|\le \frac{1}{\delta}$ ist.$\Box$

\vskip+0.3cm

{\bf 3. Die Ausnahmemenge.}

\vskip+0.3cm

F\"{u}r einen rationalen Unterraum $B$ definieren wir die Mengen
$$
\hbox{\got A} (B, \varepsilon ) =
\{  \hbox{ A} \in \hbox{\got G},\,\,\,
\psi (A,B) \le  \varepsilon\},\,\,\,\,\,
\hbox{\got A}^{[1]} (B, \varepsilon ) =
\{  \hbox{ A} \in \hbox{\got G}^{[1]}:\,\,\,
\psi (A,B) \le  \varepsilon\}.
$$
Sei $\hbox{\got Q} \subset \hbox{\got G}$ die Menge  aller rationalen zweidimensionalen Unterr\"{a}ume.
Dann ist die Ausnahmemenge von Satz 1 durch
$$
\bigcap_{H_0=1}^\infty 
\,\,\,
\bigcup_{H:\,\, H\ge H_0,\,\,
H^2\in \mathbb{Z}}
\,\,\,
\left( \bigcup_{{
 B\in \hbox{\got Q}, H(B) =H}}  \, \hbox{\got A} (B, \omega(H) )\right) 
$$
gegeben.
Sei $ \mathcal{A}_{i,j} (B,\varepsilon )  = \nu_{i,j} \left( \hbox{\got A}_{i,j}^{[1]} (B,\varepsilon)\right)$. 
Aus (\ref{00}) folgt, dass, um den Satz 1 zu beweisen, es genügt zu zeigen, dass die Mengen 
$$
\bigcap_{H_0=1}^\infty
\,\,\,
\bigcup_{H }
\,\,\, \left( \bigcup_{{
 B\in \hbox{\got Q},\, H(B) =H}}  \, \mathcal{ A}_{i,j} (B, \omega(H) )\right),\,\,\,\, i<j
$$
 Nullmengen sind.
 Mehr noch, wenn wir die Mengen
$$
\hbox{\got A} (B, \varepsilon;\delta ) =
\{  \hbox{ A} \in \hbox{\got G} (\delta):\,\,\,
\psi (A,B) \le  \varepsilon\}
$$
 betrachten, haben wir
 $$
 \bigcup_{\delta >0} \hbox{\got A} (B, \varepsilon;\delta ) = \hbox{\got A}^{[1]} (B, \varepsilon ).
 $$
  Das hei\ss t, um den Satz 1 zu beweisen genügt es zu beweisen, dass  f\"{u}r jedes $\delta >0$ die Mengen 
$$
\bigcap_{H_0=1}^\infty 
\,\,\,
\bigcup_{H }
\,\,\,
\left( \bigcup_{{
 B\in \hbox{\got Q},\, H(B) =H}}  \, \mathcal{ A}_{i,j} (B, \omega(H); \delta) \right),\,\,
 \text{wo}\,\,
 \mathcal{ A}_{i,j} (B, \varepsilon; \delta) = \mathcal{ A}_{i,j} (B, \varepsilon) \cap \mathcal{E} (\delta) $$
 Nullmengen sind.

  \vskip+0.3cm
  
   Es gilt
 $$
  \hbox{\got Q} = 
  \bigcup_{i_1<i_2}
  \hbox{\got Q}_{i_1,i_2}
  $$
  wo 
  jedes $ B \in  \hbox{\got Q}_{i_1,i_2}$
   wie in (\ref{ell})  mit 
  $
x_1 = z_{i_1}, x_2 = z_{i_2}, y_1 = z_{i_3},y_2 = z_{i_4}
$
  dargestellt werden kann, n\"{a}mlich, durch
  \begin{equation}\label{bbbbq}
  B =\{\pmb{z}= (x_1,x_2,y_1,y_2):\,\,\, y_j = b_{j,1}x_1 +
b_{j,2}x_2,\,\,\, j =1,2\},\,\,\,\,\,\,\text{mit}\,\,\,\,\,\,
\max_{i,j} |b_{i,j}|\le 1,\,\,\,\,\,\,
b_{i,j} \in \mathbb{Q}
.
\end{equation}
 Es ist klar, 
 dass $\Omega \subset  {\mathcal{E} }  (\delta)$ eine Nullmenge ist, dann und nur dann wenn
 $\nu_{i,j;k,l} (\Omega )\subset   \nu_{i,j;k,l} (\mathcal{E} (\delta)) $ eine Nullmenge ist.
    Um  Satz 1 zu beweisen, genügt es also zu zeigen, dass 
     f\"{u}r jedes $\delta >0$ 
    und f\"{u}r alle $ (k,l)$ und $(i,j)$
    die Mengen 
$$
\bigcap_{H_0=1}^\infty 
\,\,\,
\bigcup_{H}
\,\,\,
\left( \bigcup_{{
 B\in \hbox{\got Q}_{k,l},\, H(B) =H}}  \,  \nu_{i,j;k,l}  \left( \mathcal{ A}_{i,j} (B, \omega(H); \delta) \right)\right)
 $$
 Nullmengen sind.
     
Wegen $ A \in \hbox{\got G}(\delta)$
wobei die Menge $\hbox{\got G}(\delta)$  in (\ref{eee11}) definiert ist, kann $A$ als
\begin{equation}\label{aaaa}
A =
\{\pmb{z}= (x_1,x_2,y_1,y_2):\,\,\, y_j = a_{j,1}x_1 +
a_{j,2}x_2,\,\,\, j =1,2\}
\end{equation}
dargestellt  werden mit
$$
\pmb{a} = (a_{1,1},a_{1,2},a_{2,1},a_{2,2}) 
\in \Pi_\delta = 
\left\{
(a_{1,1},a_{1,2},a_{2,1},a_{2,2}) \in \mathbb{R}^4:\,\,\,\,
\max_{1\le i,j\le 2}
|a_{i,j}|
\le \frac{1}{\delta}\right\}
$$
(siehe (\ref{max})).  Dann folgt, dass
$$
 \nu_{i,j;k,l}  \left( \mathcal{ A}_{i,j} (B, \omega(H); \delta) \right)
 \subset
\Omega (B, \omega(H);  \delta)  ,
$$
wobei wir 
$$
\Omega (B,\varepsilon; \delta) 
 =
  \Pi_\delta\cap  \Omega (B,\varepsilon),\,\,
  \Omega (B,\varepsilon)  
 =\{ 
 \pmb{a}   \in \mathbb{R}^4:\,
\text{ f\"{u}r}\,
 A \,\text{aus (\ref{aaaa})  gilt }
\psi(A,B) \le \varepsilon\}
$$
gesetzt haben.
Damit haben wir bewiesen, dass
\begin{equation}\label{mera}
\mu_4\left(
\bigcap_{H_0=1}^\infty 
\,\,\,
\bigcup_{H:\,\, H\ge H_0,\,\,
H^2\in \mathbb{Z}}
\,\,\,
\left( \bigcup_{{
 B\,\,\text{aus (\ref{bbbbq})},\,\, H(B) =H}}  \,  \Omega (B, \omega(H);  \delta)
 \right)\right)=0
 \end{equation}
und Satz 1 folgt.

  \vskip+0.3cm

 {\bf 4. Die Menge   $\Omega (B, \varepsilon;  \delta)$.}
 
  \vskip+0.3cm

  {\bf Hilfssatz 3.}\,\,{\it
  Seien $ A$ und $B$ zwei Unterr\"{a}ume  definiert durch (\ref{aaaa}) und (\ref{bbbbq}).
  Dann hat mann, 
  falls $ A\cap B \neq\{ \pmb{0}\}$,
 dass
  \begin{equation}\label{deta}
  \left|
  \begin{array}{cc}
  a_{1,1} -  b_{1,1}
  &
   a_{2,1} -  b_{2,1} 
   \cr
    a_{1,2} -  b_{1,2}
  &
   a_{2,2} -  b_{2,2}
   \end{array}
  \right| = 0.
  \end{equation}
  
  }
  
    \vskip+0.3cm
  Beweis.
  Falls $ A\cap B \neq\{ \pmb{0}\}$,
  gibt es $(x_1,x_2) \neq (0,0)$ und $(y_1,y_2)$ mit
  $$
   \left(
  \begin{array}{c}
  y_1
   \cr
 y_2 
   \end{array}
  \right)
  =
  \left(
  \begin{array}{cc}
  a_{1,1}    &
   a_{2,1}  
   \cr
    a_{1,2} 
  &
   a_{2,2} 
   \end{array}
  \right)
   \left(
  \begin{array}{c}
  x_1
   \cr
 x_2 
   \end{array}
  \right)
  =
   \left(
  \begin{array}{cc}
  b_{1,1}    &
   b_{2,1}  
   \cr
    b_{1,2} 
  &
   b_{2,2} 
   \end{array}
  \right)
   \left(
  \begin{array}{c}
  x_1
   \cr
 x_2 
   \end{array}
  \right)
  .
  $$
Dann bekommt man
    $$
  \left(
  \begin{array}{cc}
  a_{1,1} -  b_{1,1}
  &
   a_{2,1} -  b_{2,1} 
   \cr
    a_{1,2} -  b_{1,2}
  &
   a_{2,2} -  b_{2,2}
   \end{array}
  \right)
   \left(
  \begin{array}{c}
  x_1
   \cr
 x_2 
   \end{array}
  \right)
  =  \left(
  \begin{array}{c}
  0
   \cr
 0
   \end{array}
  \right)
  $$
  und (\ref{deta}) folgt.$\Box$

  F\"{u}r $B$ aus (\ref{bbbb}) definieren wir die Fl\"{a}che
  \begin{equation}\label{bbbb}
  \Sigma_B = 
  \{\pmb{a}
  =
  (a_{1,1},a_{1,2},a_{2,1},a_{2,2})
  \in \mathbb{R}^4:\,\,
  \text{die Gleichung (\ref{deta}) gilt}\}.
  \end{equation}

   \vskip+0.3cm

  {\bf Hilfssatz 4.}\,\,{\it Es gilt
  $\Omega (B, \varepsilon;\delta) \subset \mathcal{U}_{\varepsilon_1} \left( \Sigma_B  \right)$
  mit
  $\varepsilon_1=
  3\varepsilon \sqrt{2+ 4\delta^{-2}}$.
  }
     \vskip+0.3cm
     
     Beweis. 
    Der Beweis verl\"{a}uft analog zu \cite{n}. 
     F\"{u}r $A$ mit $ \pmb{a}\in  \Omega (B, \varepsilon)$ gibt es $ \pmb{z}_{A} = (x_1^*,x_2^*, y_1^*, y_2^*) \in \omega_B$ mit
    $$
    \rho (A, \pmb{z}_A) \le \sigma (\omega_A, \pmb{z}_A) = \sigma(\omega_A,\omega_B)\le \varepsilon.
    $$
    Betrachten wir zwei  dreidimensionale
  Unterr\"{a}ume
  $$
  A_k = 
    \{\pmb{z}= (x_1,x_2,y_1,y_2)\in  \mathbb{R}^4_{[\pmb{z}]}:\,\,\, y_k= a_{k,1}x_1 +
a_{k,2}x_2\},\,\,\,\,\, k = 1,2,
$$
im Raum $\mathbb{R}^4 = \mathbb{R}^4_{[\pmb{z}]}$ mit Koordinaten
$\pmb{z} = (x_1,x_1,y_1,y_2)$,
so ist $A = A_1 \cap A_2$, und
  $$
    \rho (A_k , \pmb{z}_A)\le \rho (A, \pmb{z}_A) \le \varepsilon.
  $$
Daraus erh\"{a}lt man wegen $ \pmb{a} =(a_{1,2},a_{1,2}, a_{2,1},a_{2,2})\in \Pi_\delta$ die Absch\"{a}tzungen
  $$
  \frac{|y_k^* - a_{k,1}x_1^* - a_{k,2}x_2^*|}{\sqrt{1+ 2\delta^{-2}}} \le
  \frac{|y_k^* - a_{k,1}x_1^* - a_{k,2}x_2^*|}{\sqrt{1+ a_{k,1}^2+a_{k,2}^2}} \le \varepsilon,\,\,\,\,\, k = 1,2
    $$
  und
  \begin{equation}\label{lll}
  |y_k^* - a_{k,1}x_1^* - a_{k,2}x_2^*|
  \le \varepsilon \sqrt{1+ 2\delta^{-2}},\,\,\,\,\, k = 1,2.
  \end{equation}
  Nun betrachten wir dreidimensionale affine Unterr\"{a}ume
  $$
  \mathcal{Z}_k=  \mathcal{Z}_k(A)
  =
  \{
  \pmb{c} = (c_{1,1},c_{1,2},c_{2,1},c_{2,2}) \in \mathbb{R}^4_{[\pmb{c}]}:\,\,\,
  y_k^* - c_{k,1}x_1^* - c_{k,2}x_2^*
  =0\} \subset \mathbb{R}^4_{[\pmb{c}]},\,\,\,\,\, k = 1,2.
  $$
  mit den Normalvektoren
  $$
  \pmb{n}_1 = (x_1^*, x_2^*,0,0),\,\,\,\,\,\text{und}
  \,\,\,\,\,
    \pmb{n}_2 = (0,0,x_1^*, x_2^*),
  $$
  die orthogonal zueinander sind. 
  Diese Unterr\"{a}ume sind im euklidischen  Raum $\mathbb{R}^4 = \mathbb{R}^4_{[\pmb{c}]}$ mit Koordinaten $(c_{1,1},c_{1,2},c_{2,1},c_{2,2})$ enthalten.
  Man hat
 $$
  \rho (\pmb{a}, \mathcal{Z}_k) =
  \frac{  |y_k^* - a_{k,1}x_1^* - a_{k,2}x_2^*|}{\sqrt{|x_1^*|^2+|x_2^*|^2}},\,\,\,\,\, k = 1,2
  .
  $$
Da
  $
  \pmb{z}_A \in \omega_B \subset \hbox{\got S},$
 $
 |b_{i,j}|\le 1,
 $
 haben wir
 $$
 1=
 {|x_1^*|^2+|x_2^*|^2+
 |y_1^*|^2+|y_2^*|^2
  }
  \le
 |x_1^*|^2+|x_2^*|^2 + 8\max_{j=1,2}|x_j|^2\le
 9 ( |x_1^*|^2+|x_2^*|^2)
  ,
  $$
   und
    $$
  \rho (\pmb{a}, \mathcal{Z}_k) \le
  3 |y_k^* - a_{k,1}x_1^* - a_{k,2}x_2^*| 
  \le 3\varepsilon \sqrt{1+ 2\delta^{-2}}
  ,\,\,\,\,\, k = 1,2
  .
  $$
  Nun gilt  f\"{u}r die Menge
  $\mathcal{Z} =\mathcal{Z}(A)
  = \mathcal{Z}_1\cap 
  \mathcal{Z}_2$,
  dass
  $$
  \rho (\pmb{a}, \mathcal{Z}(A)) =
  \sqrt{ \rho (\pmb{a}, \mathcal{Z}_1)^2+  \rho (\pmb{a}, \mathcal{Z}_2)^2 }\le
  3\varepsilon \sqrt{2+ 4\delta^{-2}} =\varepsilon_1.
 $$
 Es gibt daher ein $ \pmb{c}^* \in \mathcal{Z}(A)$ mit
 $$
 \rho (\pmb{a}, \pmb{c}^*) =
 \rho (\pmb{a}, \mathcal{Z}(A)) \le
 \varepsilon_1.
 $$
 F\"{u}r jeden Unterraum $ C\subset  \mathbb{R}^4_{[\pmb{z}]}$ mit $\pmb{c} \in \mathcal{Z} =
   \mathcal{Z}(A)$ haben wir
 $ \pmb{z}_A \in C \cap B$, sodass
 $C \cap B \neq\{\pmb{0}\}$ und
  $ \pmb{c} \in \Sigma_B$ nach dem Hilfssatz 3 gilt.
 Also  ist $ \pmb{c}^* \in \Sigma_B$ und $\rho (\pmb{a}, \Sigma_B) \le 
 \rho (\pmb{a},\pmb{c}^*)\le\varepsilon_1$.$\Box$
 
    \vskip+0.3cm
    
    Es ist klar, dass
 $$
 \Sigma_B = \Sigma_0 +\pmb{b},
 \,\,\,\,
 $$
 wobei
 $$
 \pmb{b}
 =
 (b_{1,1},b_{1,2},b_{2,1},b_{2,2}) \in \mathbb{R}^4
 $$
 und
 $$
 \Sigma_0 = \{
  (c_{1,1},c_{1,2},c_{2,1},c_{2,2}) 
 \in \mathbb{R}^4:\,\,\,
 c_{1,1} c_{2,2} - c_{1,2} c_{21} =0\}.
 $$
 Es erweist sich als praktisch, mit anderen Koordinaten  
 $$
 \xi_1 = \frac{c_{1,1}+c_{2,2}}{2},\,\,\,\,\,
  \eta_1 = \frac{c_{1,1}-c_{2,2}}{2},\,\,\,\,\,
   \xi_2=  \frac{c_{2,1}-c_{1,2}}{2},\,\,\,\,\,
     \eta_2 = 
     \frac{c_{2,1}+c_{1,2}}{2}
 $$
 zu arbeiten.
 Dann gilt
 $$
 c_{1,1} c_{2,2} = \xi_1^2-\eta_1^2,\,\,\,\,\,
  c_{1,2} c_{2,2} =  \eta_2^2- \xi_2^2,
 $$
 und  statt der 
 Fl\"{a}che
 $
 \Sigma_0 $
   k\"{o}nnen  wir die  Fl\"{a}che
 $$
 \Sigma
 =\{ (\xi_1,\xi_2,\eta_1,\eta_2)\in \mathbb{R}^4:\,\,
 \xi_1^2 + \xi_2^2 = \eta_1^2+ \eta_2^2\}
 $$
 betrachten.

       \vskip+0.3cm
    {\bf 5.  \"Uber die Umgebung der  Fl\"{a}che  $
 \Sigma  
 $.}
    \vskip+0.3cm

 Betrachten wir nun die Kugel 
 $$
 \mathcal{ K}_T = \{\pmb{\zeta}= (\xi_1,\xi_2,\eta_1,\eta_2) \in \mathbb{R}^4:\,\,\,
 \xi_1^2+
 \xi_2^2+\eta_1^2+
 \eta_2^2 \le T^2\}.
 $$
 Damit ergibt sich folgender Hilfssatz.

   \vskip+0.3cm
 {\bf Hilfssatz 5.}\,{\it 
 F\"ur jede $\varepsilon >0$ , $ T \ge 1$ und $ \pmb{Z} \in \mathbb{R}^4$ gilt
 
 $$
   \mu_4\left((\mathcal{ U}_\varepsilon (\pmb{Z}+\Sigma) )\cap \mathcal{ K}_T \right)  \ll T^3\varepsilon +\varepsilon^4.
 $$
 }
    \vskip+0.3cm
 
 Beweis.   
 \,
 Wir
 betrachten die Parameterdarstellung
 $$
 \pmb{\zeta} (r,\varphi,\psi)
 :\,\,\,\,
 \xi_1 = r\cos \varphi,\,\,\, \xi_2 = r\sin\varphi,\,\,\,\, \eta_1 = r\cos \psi,\,\,\, \eta_2 = r\sin\psi
 $$
 der  Fl\"{a}che $\Sigma$.
 Dann ist 
 $$
 (\pmb{Z}+ \Sigma )\cap  \mathcal{ K}_T = \{ \pmb{\zeta } =
 \pmb{Z}+\pmb{\zeta} (r,\varphi,\psi),\,\,\,
  (r,\varphi,\psi)\in \mathcal{D}_{T} \}
 $$ 
 das  Bild eines  beschr\"{a}nkten Gebietes $\mathcal{D}_T=  \mathcal{D}_{T, \pmb{Z}}  \subset \mathbb{R}^3$.
 Betrachten wir 
 den Normalenvektor
 $$
 \pmb{n} =\pmb{n}_{\pmb{\zeta}}=
  \pmb{n} (r,\varphi,\psi) =
  \left( 
  \begin{array}{c}
  n_1\cr
  n_2\cr
  n_3
  \cr n_4
  \end{array}
  \right)
  =
   \left( 
  \begin{array}{c}
 - \frac{\xi_1}{\sqrt{2}r}\cr
-\frac{\xi_2}{\sqrt{2}r}\cr
\frac{\eta_1}{\sqrt{2}r}\,
  \cr
  \frac{\eta_2}{\sqrt{2}r}
  \end{array}
  \right)=
    \left( 
  \begin{array}{c}
  -\frac{\cos \varphi}{\sqrt{2}}\cr
-\frac{\sin \varphi }{\sqrt{2}}\cr
\frac{\cos \psi }{\sqrt{2}}\,
  \cr
  \frac{\sin \psi }{\sqrt{2}}
  \end{array}
  \right),\,\,\,\,\,
  |\pmb{n}| = 1
  $$
  f\"{u}r $\pmb{Z}+\Sigma $ im Punkt
  $\pmb{\zeta } =\pmb{Z}+\pmb{\zeta} (r,\varphi,\psi)$.
 
 Die Funktion 
 $$
 \hbox{\got f}:\,\,
  \mathbb{R}_+\times (\mathbb{R}/2\pi \mathbb{Z})^2\times [-\varepsilon,+\varepsilon] \to \mathcal{ U}_\varepsilon (\pmb{Z}+\Sigma) ,\,\,\,\,\,\,
 \pmb{q} = (r,\varphi,\psi, t) \mapsto \pmb{\zeta}_t = \pmb{\zeta}+ t \pmb{n}_{\pmb{\zeta}}
 $$
 gibt eine Parameterdarstellung von $ \mathcal{ U}_\varepsilon (\pmb{Z}+\Sigma) $. Hier wollen wir bemerken, dass
 diese Funktion  surjektiv, aber nicht injektiv ist,   und so  hat man f\"{u}r $\mathcal{D}\subset \mathbb{R}^3$
 dass
 \begin{equation}\label{a}
 \mu_4 \left(\mathcal{F} (\mathcal{D} \times [-\varepsilon,+\varepsilon])\right)
 \le
 \int_{\mathcal{D}} \int_{-\varepsilon}^{+\varepsilon}
 \sqrt{{\rm det}\,G} \, dr\, d\varphi\, d\psi \, dt
 ,\end{equation}
 wobei
 $G$ ist die Gramsche Matrix der Vektoren
 $$
 \frac{\partial  \pmb{\zeta}_t }{\partial t}
 = \pmb{n}_{\pmb{\zeta}} 
 ,\,\,\,
  \frac{\partial \pmb{\zeta}_t }{\partial r}
  =
  \left( 
  \begin{array}{c}
 {\cos \varphi} 
 \cr
  \sin\varphi  
  \cr
   \cos\psi 
  \cr 
   \sin\psi  
  \end{array}
  \right)
  ,\,\,\,
  \frac{\partial  \pmb{\zeta}_t }{\partial \varphi}
  =
  \left( 
  \begin{array}{c}
  \left( \frac{t}{\sqrt{2}}- r\right) \sin\varphi
  \cr
   \left( r-\frac{t}{\sqrt{2}}\right) \cos\varphi
   \cr
   0
   \cr
   0
 \end{array}
  \right)
  ,\,\,\,
  \frac{\partial  \pmb{\zeta}_t }{\partial \psi}=
  \left( 
  \begin{array}{c}
  0
  \cr
  0
  \cr
  -\left( r+\frac{t}{\sqrt{2}}\right) \sin\psi
  \cr
   \left( r+\frac{t}{\sqrt{2}}\right) \cos\psi
  \end{array}
  \right).
 $$ 
 Nun erh\"{a}lt man
 $$
 G
 =
 \left|
 \begin{array}{cccc}
 1 &  0 & 0  & 0   \cr
  0  & 2  &  0 & 0   \cr
   0 &   0 & \left( r-\frac{t}{\sqrt{2}}\right)^2&  0  \cr
   0 &    0 &    0       & \left( r+\frac{t}{\sqrt{2}}\right)^2
 \end{array}
 \right| = 2 \left(r^2-\frac{t^2}{2}\right)^2,
 $$
 und f\"{u}r $|t| \le \varepsilon, r\ge  {\varepsilon} $
 folgt daraus, dass
 \begin{equation}\label{c}
 |G|\le r^4.
 \end{equation} 
  Aus der Dreiecksungleichung folgt nun, dass 
 $$
 \mathcal{ U}_\varepsilon (\pmb{Z}+\Sigma) \cap \mathcal{ K}_T \subset
 \hbox{\got f} \left(\mathcal{D}_{T+\varepsilon} \times [-\varepsilon,+\varepsilon]\right)
  $$
  gilt.
 Tats\"{a}chlich, 
 falls $\pmb{\gamma} \in \mathcal{ U}_\varepsilon (\pmb{Z}+\Sigma) \cap \mathcal{K}_T$ ist, so  gilt
 f\"{u}r den Punkt $\pmb{\zeta} \in  \pmb{Z}+\Sigma  $  und $ t \in [-\varepsilon,+\varepsilon]$ mit $ \pmb{\gamma} = \pmb{\zeta}+ t\pmb{n}_{\pmb\zeta} \in \mathcal{ K}_T$ 
 die Ungleichung
 $ |\pmb{\zeta}-\pmb{\gamma}|\le \varepsilon$. Also haben wir
 $\pmb{\gamma} \in \mathcal{K}_{T+\varepsilon}$.
 
 Wir betrachten nun die zwei Mengen
 $$
 \mathcal{D}_T^{[1,\varepsilon]} =
 \{ (r,\varphi,\psi) \in \mathcal{D}_T:\,\,\, r\ge {\varepsilon} \}\,\,\,\,\,
  \text{und}\,\,\,\,
   \mathcal{D}_T^{[2,\varepsilon]}=  \mathcal{D}_T\setminus \mathcal{D}_T^{[1,\varepsilon]} 
   $$
   und ihre Bilder
   $$
   \mathcal{W}^{[j,\varepsilon]}_T=
   \hbox{\got f} \left(\mathcal{D}_T^{[j,\varepsilon]} \times [-\varepsilon,+\varepsilon])\right),\,\,\,\,\, j =1,2.
  $$
  F\"{u}r sie  haben wir  
  \begin{equation}\label{1q}
   \mu_4\left( \mathcal{ U}_\varepsilon (\pmb{Z}+\Sigma) \cap \mathcal{ K}_T 
   \right)\le
\mu_4\left(\hbox{\got f} (\mathcal{D}_{T+\varepsilon} \times [-\varepsilon,+\varepsilon])
  \right)
  =
  \mu_4\left(\mathcal{ W}_{T+\varepsilon}^{[1,\varepsilon]} \right)
  +
 \mu_4\left(\mathcal{ W}_{T+\varepsilon}^{[2,\varepsilon]} \right)
  .
  \end{equation}
  F\"{u}r $\mathcal{D} = \mathcal{D}_{T+\varepsilon}^{[1,\varepsilon]}$ aus (\ref{a}) und (\ref{c}) haben wir
  die Ungleichung
  \begin{equation}\label{2q}
  \mu_4\left(\mathcal{ W}_{T+\varepsilon}^{[1,\varepsilon]} \right)\le
 \int_{\mathcal{D}_{T+\varepsilon}^{[1,\varepsilon]}} \int_{-\varepsilon}^{+\varepsilon}
 r^2 \, dr\, d\varphi\, d\psi \, dt\ll T^3\varepsilon.
 \end{equation}
  
  Nun sch\"{a}tzen wir die Werte von
  $\mu_4\left(\mathcal{ W}_{T+\varepsilon}^{[2,\varepsilon]} \right)
$ ab.
 Aus der Dreiecksungleichung  folgt
 $$
 \mathcal{W}_{T+\varepsilon}^{[2,\varepsilon]}=
 \hbox{\got f} (\mathcal{D}_{T+\varepsilon}^{[2,\varepsilon]}\times[-\varepsilon,+\varepsilon]) \subset 
 \mathcal{ U}_{ 2\varepsilon}(\pmb{Z}),
 $$
 sodass
 \begin{equation}\label{3q}
 \mu_4\left(\mathcal{ W}_{T+\varepsilon}^{[2,\varepsilon]} \right)
 \le
 \mu_4\left(
  \mathcal{U}_{ 2{\varepsilon}}(\pmb{Z})
 \right)
 \ll \varepsilon^4.
 \end{equation}
 Hilfssatz folgt  nun aus (\ref{1q},\ref{2q},\ref{3q}).$\Box$
 
    \vskip+0.3cm
 
 {\bf 6. Beweis von Satz 1.}
 
    \vskip+0.3cm
 
 Es  gen\"{u}gt, (\ref{mera})  zu beweisen.
 Aus Hifssatz 4, Hilfssatz 5 mit $ T \asymp \delta^{-1}, \varepsilon  \asymp \omega(H) \sqrt{1+2\delta^{-2}}
 \asymp \omega(H)$
 und  der Definition der Fl\"{a}che $\Sigma$   folgt, dass, 
 um (\ref{mera})  zu beweisen, 
 man zeigen muss, dass die Reihe
 \begin{equation}\label{reiche}
 \sum_{H: \, H\ge 1,\, H^2\in \mathbb{Z}}\,\,\,
 \sum_{{
 B\,\,\text{aus (\ref{bbbbq})},\,\, H(B) =H}}
 \omega(H).
 \end{equation}
 konvergiert. Sei dazu
 $$
 N(H)=
 \#\{ B \,\text{rationaler zwei-dimensionaler Unterraum in }\, \mathbb{R}^4:\,\,
 H(B) = H\}.
 $$
 Es ist bekanntlich so, dass
 $$
 \sum_{1\le k\le H, \, k^2\in \mathbb{Z}} \,\,\,N(k) =
  \sum_{1\le l\le H^2, \, l\in \mathbb{Z}} \,\,\,N\left(\sqrt{l}\right) =
 $$
 $$
 =
 \#\{ B \,\text{rationaler zwei-dimensionaler Unterraum in }\, \mathbb{R}^4:\,\,
 H(B) \le H\} \asymp H^4.
 $$
 (siehe Theorem 3 aus \cite{s1}).
 Mit Hilfe der 
  abelschen partiellen Summation erhalten wir
 $$
  \sum_{
  1\le H \le \sqrt{W},\, H^2 \in \mathbb{Z}
  }\,\,\,
 \sum_{{
 B\,\,\text{aus (\ref{bbbbq})},\,\, H(B) =H}}
  \omega(H)
 =
 \sum_{j=1}^W
 N\left(\sqrt{j}\right)
  \omega\left(\sqrt{j}\right)=
 $$
 $$=
  \sum_{j=1}^{W-1}\,
  \left(\sum_{l\le j}
 N\left(\sqrt{l}\right) \right)\left(\omega\left(\sqrt{j}\right)- \omega\left(\sqrt{j+1}\right)\right)+
 \omega\left(\sqrt{W}\right) 
  \sum_{l\le W}
 N\left(\sqrt{l}\right)\ll
 $$
 $$
 \ll
   \sum_{j=1}^{W-1}\,
  j^2
  \left(\omega\left(\sqrt{j}\right)- \omega\left(\sqrt{j+1}\right)\right)
  +
   \omega\left(\sqrt{W}\right) 
 W^2  
 \ll
 \sum_{j=1}^W  j \, \omega \left(\sqrt{j}\right).
 $$
 Daraus folgt, dass
  (\ref{reiche}) konvergiert, falls  (\ref{reiche0}) konvergiert.$\Box$

     \vskip+0.3cm
     
     {\bf 7. Eine Verallgemeinerung.}
 
       \vskip+0.3cm
       
       Hier möchten wir eine allgemeine Aussage formulieren, die mit der Methode dieses Artikels bewiesen werden kann.
       Wir formulieren sie ohne Beweis.
        \vskip+0.3cm
  
    {\bf Satz 2.}\,\,{\it  Sei $d=2s$. 
    Sei $\omega (j) \ge 0, j \ge 1$
    eine monoton fallende Funktion. Angenommen die Reihe
    $$
 \sum_{H=1}^\infty 
    j^{s-1}\omega \left(\sqrt{j}\right)
$$
 konvergiert.
   Dann gibt es  f\"{u}r fast jeden  $s$-dimensionalen Unterraum $A\subset\mathbb{R}^d$ 
    eine positive Konstante $C(A)$ mit
   $$
     \psi (A,B) \ge C(A) \, \omega (H(B)) \,\,\,\,\,
    \text{f\"{u}r alle $s$-dimensionalen rationalen Unterr\"{a}ume}\,\,\,
    B
    \,\,\,\text{in}\,\,\,
    \mathbb{R}^d
    .
$$
    }
 
       \vskip+0.3cm

 {\bf Bemerkung. }     F\"{u}r $n=k=s,d= 2s$ folgt aus Satz B, dass  es einen $s$-dimensionalen Unterraum $A$ mit
    $$
     \psi (A,B) \ge C H^{-s^2-1},\,\,\,\,\,
    \text{f\"{u}r alle $s$-dimensionalen rationalen Unterr\"{a}ume}\,\,\,
    B
      \,\,\,\text{in}\,\,\,
    \mathbb{R}^d
    $$
    gibt. Aus Satz 2 folgt 
    die Versch\"{a}rfung: f\"{u}r jedes $\varepsilon >0$
   gibt es einen $s$-dimensionalen Unterraum $A$ mit
      $$
       \psi (A,B) \ge C(A,\varepsilon) 
        H^{-2s-\varepsilon},\,\,\,\,\,
    \text{f\"{u}r alle $s$-dimensionalen rationalen Unterr\"{a}ume}\,\,\,
    B   \,\,\,\text{in}\,\,\,
    \mathbb{R}^d.
      $$

        \vskip+0.3cm 
       
       Der Beweis dieses allgemeines Resultates verl\"{a}uft analog zu jenem von Satz 1.
       Die Hauptunterschied besteht darin, dass statt der quadratischen Fläche (\ref{bbbb})   man f\"{u}r die Matrix
       $$
       B = 
       \left(
       \begin{array}{ccc}
       b_{1,1}& ... & b_{1,s}\cr
       ...& ... & ... \cr
       b_{s,1}& ... & b_{s,s}
       \end{array}
       \right)
       $$
       eine nicht quadratische  Fl\"{a}che
       $$
        \Sigma_B^s = 
  \left\{ A =
         \left(
       \begin{array}{ccc}
       a_{1,1}& ... & a_{1,s}\cr
       ...& ... & ... \cr
       a_{s,1}& ... & a_{s,s}
       \end{array}
       \right)
  \in \mathbb{R}^{s^2}:\,\,\,\,
  {\rm det}\, ( A- B) = 0
  \right\} \subset \mathbb{R}^{s^2}
 $$
 betrachten  und die Ungleichung
 $$
 \mu_{s^2} \left( \mathcal{U}_\varepsilon \left(
 \Sigma_B^s \cap
 [-T,T]^{2s}\right)\right) \ll_T \varepsilon
 $$
 beweisen muss.

 \vskip+1cm
 
 Autor: Nikolay Moshchevitin,
 
 Nationale Pazifik-Universit\"{a}t, Russland;
 
 e-mail: moshchevitin@rambler.ru, moshchevitin@gmail.com

 \end{document}